\documentclass[11pt]{amsart}
\textheight
9in
\usepackage{amsmath,amsthm,amsfonts,amssymb,latexsym,epsfig}

\newtheorem{theorem}{Theorem}[section]

\newtheorem{lemma}{Lemma}[section]
\newtheorem{prop}{Proposition}[section]

\newtheorem{cor}{Corollary}[section]

\numberwithin{equation}{section}

\theoremstyle{definition}

\theoremstyle{remark}

\begin{document}
\title{Some Monotonicity Properties of Gamma and $q$-gamma Functions}
\author{Peng Gao}
\address{Department of Computer and Mathematical Sciences,
University of Toronto at Scarborough, 1265 Military Trail, Toronto
Ontario, Canada M1C 1A4}
\email{penggao@utsc.utoronto.ca}
\date{September 21, 2007.}

\begin{abstract}
We prove some properties of completely monotonic
functions and apply them to obtain results on gamma and
$q$-gamma functions.
\end{abstract}

\maketitle
\section{Introduction}
\label{sec 1} \setcounter{equation}{0}
   The $q$-gamma function is defined for positive real numbers $x$ and $q \neq 1$ by
\begin{eqnarray}
\label{01}
  \Gamma_q(x)&=&
      (1-q)^{1-x}\prod^{\infty}_{n=0}\frac {1-q^{n+1}}{1-q^{n+x}}, ~~
      0<q<1; \\
  \Gamma_q(x)&=&(q-1)^{1-x}q^{\frac 1{2}x(x-1)}\prod^{\infty}_{n=0}\frac {1-q^{-(n+1)}}{1-q^{-(n+x)}}
  , ~~q>1. \nonumber
\end{eqnarray}
    We note here \cite{Koo} the limit of $\Gamma_q(x)$ as $q \rightarrow
    1^{-}$ gives back the well-known Euler's gamma function:
\begin{equation*}
    \lim_{q \rightarrow
    1^{-}}\Gamma_q(x)=\Gamma(x)=\int^{\infty}_0 t^x e^{-t}\frac
    {dt}{t}.
\end{equation*}
    Note also from the
    definition, for positive $x$ and $0<q<1$,
    $$\Gamma_{1/q}(x)=q^{(x-1)(1-x/2)}\Gamma_{q}(x),$$
we see that $\lim_{q \rightarrow 1}\Gamma_q(x)=\Gamma(x)$. For
historical remarks on gamma and $q$-gamma functions, we
    refer the reader to \cite{Koo}, \cite{alz1} and \cite{alz1.5}.

    There exists an extensive and rich literature on inequalities for the gamma and $q$-gamma functions.
For the recent developments in this area, we refer the reader to the
articles \cite{I&M}, \cite{alz1}, \cite{alz1.5}-\cite{alz2}, \cite{Q&V} and the
references therein. Many of these inequalities follow from the
monotonicity properties of functions which are closely related to
$\Gamma$ (resp. $\Gamma_q$) and its logarithmic derivative $\psi$
(resp. $\psi_q$). Here we recall that a function $f(x)$ is said to
be completely monotonic on $(a, b)$ if it has derivatives of all
orders and $(-1)^kf^{(k)}(x) \geq 0, x \in (a, b), k \geq 0$. Lemma
\ref{lem1} below asserts that $f(x)=e^{-h(x)}$ is completely
monotonic on an interval if $h'$ is. Following \cite{G&I}, we call
such functions $f(x)$ logarithmically completely monotonic.

    We note here that $\lim_{q \rightarrow 1}\psi_q(x)=\psi(x)$ (see
     \cite{K&S})
     and that $\psi'$ and $\psi'_q$ are completely
    monotonic functions on $(0, +\infty)$ (see \cite{K&M}, \cite{alz2}). Thus, one expects to deduce results on
    gamma and $q$-gamma functions from properties of
(logarithmically) completely monotonic functions,
    by applying them to functions related to $\psi'$ or $\psi'_q$.
    It is our goal in this paper to obtain some results on gamma and
$q$-gamma functions via this approach.
   As an example, we
recall the following result of Bustoz and Ismail (the case $q=1$) as
well as Ismail and Muldoon:
\begin{theorem}[{\cite[Theorem 3]{B&I}, \cite[Theorem 2.5]{I&M}}]
\label{thm5}
   Let $a+1 \geq b >a$, $\alpha=\max (-a, -c)$
   and for $q>0$, define
\begin{equation*}
   g_q(x;a,b,c)=\Big(\frac {1-q^{x+c}}{1-q} \Big )^{a-b}\frac {\Gamma_q(x+b)}{\Gamma_q(x+a)},
   \hspace{0.1in} x>\alpha,
\end{equation*}
  where $\Gamma_1(x)=\Gamma(x)$. Then $g_q(x;a,b,c)$ is logarithmically completely monotonic on $(\alpha,
  +\infty)$ if $c \leq (a+b-1)/2$ and $1/g_q(x;a,b,c)$ is logarithmically completely monotonic
   on $(\alpha,
  +\infty)$ if $c \geq a$.
\end{theorem}
   It follows immediately from the above theorem that for $0< q< 1$ and
$0 < s <1$, one has \cite{I&M} for $x>0$:
\begin{equation*}
    \Big(\frac {1-q^{x+s/2}}{1-q} \Big )^{1-s} < \frac
    {\Gamma_q(x+1)}{\Gamma_q(x+s)} < \Big(\frac {1-q^{x+s}}{1-q} \Big
    )^{1-s}.
\end{equation*}
   Alzer later \cite{alz1.5} determined the best values $u(q,s), v(q,s)$ such that the inequalities:
\begin{equation}
\label{1.4}
    \Big(\frac {1-q^{x+u(q,s)}}{1-q} \Big )^{1-s} < \frac
    {\Gamma_q(x+1)}{\Gamma_q(x+s)} < \Big(\frac {1-q^{x+v(q,s)}}{1-q} \Big
    )^{1-s}
\end{equation}
   hold for $q>0, 0<s<1, x>0$ to be
\begin{eqnarray}
\label{1.4'}
   u(q, s) &=& \frac {\ln \frac {q^s-q}{(1-s)(1-q)}}{\log q},
\hspace{0.2in} 0<q<1;  \\
   \alpha(q, s) &=& s/2,
\hspace{0.2in} 0<q<1; \nonumber \\
   v(q, s) &=& \frac {\ln \Big (1-(1-q)\Gamma^{1/(s-1)}_q(s) \Big )}{\ln q}. \nonumber
\end{eqnarray}

   Motivated by the above results, we will show that the function $g_q(x; s, 1, u(q, s))$ as defined in
   Theorem \ref{thm5} is
logarithmically completely monotonic on $(0,
  +\infty)$ for $0<q<1$. This will enable us to deduce the left-hand
  side inequality of Alzer's result above.

   Let $P_{n,k} \hspace{0.02in} (1 \leq k \leq n)$ be the set of all
vectors ${\bf
   m}=(m_1, \ldots, m_k)$ whose components are natural numbers such
   that $1\leq m_{\nu} < m_{\mu} \leq n$ for $1 \leq \nu < \mu \leq
   k$ and let $P_{n,0}$ be the empty set. Recently,  Grinshpan and Ismail \cite{G&I} proved
   the following result:
\begin{theorem}[{\cite[Theorem 1.2, 3.3]{G&I}}]
\label{thm3.0}
   For any $a_k>0 \hspace{0.02in} (1 \leq k \leq n)$, define for
   $x>0, 0< q \leq 1$,
\begin{equation*}
   F_n(x, q)=\frac {\Gamma_q(x)\prod^{[n/2]}_{k=1}\Big (\prod_{{\bf m} \in
   P_{n,2k}}\Gamma_q(x+\sum^{2k}_{j=1}a_{m_j}) \Big
)}{\prod^{[(n+1)/2]}_{k=1}\Big (\prod_{{\bf m} \in
   P_{n,2k-1}}
   \Gamma_q(x+\sum^{2k-1}_{j=1}a_{m_j})  \Big )}.
\end{equation*}
  Then $F_n(x, q)$ is a logarithmically completely monotonic function of $x$ on $(0,
  +\infty)$. Here we define $\Gamma_1(x)=\Gamma(x)$.
\end{theorem}
   We will show in Section \ref{sec 4.1} that the above theorem follows from a simple
   observation on logarithmically completely monotonic functions which also
enables us to deduce some other useful inequalities.

   The derivatives of $\psi(x)$ are known as polygamma functions and
   in Section \ref{sec 4.2} we will prove some inequalities
   involving the polygamma functions.

    Formulas for the volumes of geometric bodies sometimes involve the gamma
function. One then expects to obtain inequalities involving such
volumes via relevant inequalities involving gamma functions. In
Section \ref{sec 5}, we will have a brief discussion on inequalities
for the volume of the unit ball in $R^n$.

\section{Lemmas}
\label{sec 3} \setcounter{equation}{0}
   We will need the following two facts about $\Gamma(x)$ and $\psi_q(x)$. These can be found, for example, in
  \cite[(7.1)]{alz1}, \cite[(2.7)]{alz1.5} and \cite[(1.2)-(1.5)]{alz2}:
\begin{lemma}
\label{lem3.0}
  For $x>0$,
\begin{eqnarray}
\label{3.05}
   \ln \Gamma(x) &=& (x-\frac 1{2})\ln x-x+\frac 1{2}\ln(2\pi)+O\genfrac(){1pt}{}
  {1}{x}, \hspace{0.1in} x \rightarrow +\infty, \\
\label{1.12}
   \psi_q(x) & = &  -\ln (1-q) + \ln q \sum^{\infty}_{n=1}\frac {q^{nx}}{1-q^n},  \hspace{0.1in}
   0<q<1, \\
\label{1.10}
  (-1)^{n+1}\psi^{(n)}(x) &=& \int^{\infty}_0e^{-xt}\frac {t^n}{1-e^{-t}}dt
  =n!\sum^{\infty}_{k=0}\frac 1{(x+k)^{n+1}}, \hspace{0.1in} n \geq
  1, \\
\label{3.4}
  \psi^{(n)}(x+1) &=& \psi^{(n)}(x)+(-1)^n\frac {n!}{x^{n+1}},
  \hspace{0.1in} n \geq 0, \\
\label{3.03}
  \psi(x) &=& \ln x - \frac 1{2x}-\frac 1{12x^2}+O\genfrac(){1pt}{}{1}{x^3}, \hspace{0.1in} x \rightarrow +\infty,\\
\label{1.2}
  (-1)^{n+1}\psi^{(n)}(x) &=& \frac {(n-1)!}{x^n}+\frac{n!}{2x^{n+1}}+O\genfrac(){1pt}{}
  {1}{x^{n+2}}, \hspace{0.1in} n \geq 1, \hspace{0.1in} x \rightarrow +\infty.
\end{eqnarray}
\end{lemma}

\begin{lemma}[{\cite[Lemma 2.2]{Gao}}]
\label{thm11} For fixed $n \geq 1, a \geq 0$, the function
$f_{a,n}(x)=x^n(-1)^{n+1}\psi^{(n)}(x+a)$ is increasing on $[0,
+\infty)$
  if and only if $a \geq 1/2$. Also, $f_{0,n}(x)$ is decreasing on $(0, +\infty)$.
\end{lemma}
   We note here we may regard Lemma \ref{thm11} as saying that both
   $x\psi'(x)$ and $\psi'(x+a)+x\psi''(x+a), a \geq 1/2$ are completely monotonic on
   $(0, +\infty)$.

    The set of all the completely monotonic functions on an interval can be
  equipped with a ring structure with the usual addition and multiplication of functions. The next two simple lemmas
  can also be used to construct new (logarithmically) completely monotonic function, part of Lemma \ref{thm6} is contained in \cite[Lemma 1.3]{I&M}:
\begin{lemma}
\label{thm6}
   If $f(x)$ is completely monotonic on some interval $(a, b)$, then so is $f(x)-f(x+c)$ on
   $(a,b) \cap (a-c, b-c)$ for
   any $c>0$. Consequently, if $f(x)$ is logarithmically completely monotonic on some interval $(a, b)$, then so is $f(x)/f(x+c)$ on
   $(a,b) \cap (a-c, b-c)$ for
   any $c>0$.
\end{lemma}
\begin{lemma}[{\cite[Lemma 2.1]{B&I}}]
\label{lem1}
   If $f'(x)$ is completely monotonic on an interval, then
$\exp(-f(x))$ is also completely monotonic on the same interval.
\end{lemma}

\begin{lemma}
\label{lem2} Let $a_i$ and $b_i$ $(i=1,\dotsc,n)$ be real numbers
such that $0<a_1\le\dotsb\le a_n$, $0<b_1\le\dotsb\le b_n$, and
$\sum_{i=1}^k a_i\le\sum_{i=1}^k b_i$ for $k=1,\dotsc,n$. If the
function $f(x)$ is decreasing and convex on $(0, +\infty)$, then
\[\sum_{i=1}^n f(b_i)\le\sum_{i=1}^n f(a_i).\]
  If $\sum_{i=1}^n a_i=\sum_{i=1}^n b_i$, then one only needs
$f(x)$ to be convex
  for the above inequality to hold.
\end{lemma}
  The above lemma is similar to Lemma 2 in \cite{alz1}, except here
  we assume $a_i, b_i$'s to be positive and $f(x)$ defined on
  $(0, +\infty)$. We leave the proof to the reader by pointing out that
  it follows from the theory
  of majorization, for example, see the discussions in Chap. 1, $\S 28 -\S 30$
  of \cite{B&B}.

\begin{lemma}[Hadamard's inequality]
\label{lem3}
   Let $f(x)$ be a convex function on $[a, b]$, then
\begin{equation*}
   f(\frac {a+b}2) \leq \frac {1}{b-a}\int^b_a f(x)dx \leq \frac {f(a)+f(b)}{2}.
\end{equation*}
\end{lemma}

\begin{lemma}[{\cite[Lemma 2.1]{alz1.5}}]
\label{lem4}
   Let $a > 0, b > 0$ and $r$ be real numbers with $a \neq b$, and
   let
\begin{eqnarray*}
   L_r (a, b) &=& \genfrac(){1pt}{}{a^r - b^r}{r(a-b)}^{1/(r-1)}  \hspace{0.2in} (r \neq  0,
   1), \\
   L_0(a,b) &=& \frac {a - b}{ \log a-\log
  b}, \\
   L_1(a,b) &=& \frac 1{e}\genfrac(){1pt}{}{a^a}{b^b}^{1/(a-b)}.
\end{eqnarray*}
   The function $r \mapsto L_r(a,b)$ is strictly increasing on ${\mathbb R}$.
\end{lemma}

\begin{lemma}
\label{lem10}
   For real numbers $0<s, q<1$ and an integer $n \geq 1$,
\begin{equation*}
  \Big (\frac {q^s-q}{(1-s)(1-q)} \Big )^n \geq \frac {q^{ns}-q^n}{(1-s)(1-q^n)}.
\end{equation*}
\end{lemma}
\begin{proof}
  We define for $x \geq 1$,
\begin{equation*}
  f(x)=x\ln \frac {q^{s-1}-1}{(1-s)(1-q)}  - \ln \frac
  {q^{x(s-1)}-1}{(1-s)(1-q^x)}.
\end{equation*}
  It suffices to show that $f'(x) \geq 0$ as $f(1)=0$. Now
\begin{equation*}
  f'(x)= \ln \frac {q^{s-1}-1}{(1-s)(1-q)} - \ln q^{s-1}- \frac {\ln
  q^{s-1}}{q^{x(s-1)}-1}+\ln q -\frac {\ln q}{1-q^x}.
\end{equation*}
  Note that
\begin{equation*}
  f''(x)= \frac { (\ln q)^2
  (1-s)^2q^{x(1-s)}}{(1-q^{x(1-s)})^2} -\frac {(\ln q)^2 q^x}{(1-q^x)^2}.
\end{equation*}
   If we set $u=1-s, y=q^x$ so that $0 < u < 1$ and $0 < y
   <1$ and set
\begin{equation*}
   g(u)=\frac {u^2y^u}{(1-y^u)^2}.
\end{equation*}
   Then we have
\begin{equation*}
   g'(u)= \frac {uy^u}{(1-y^u)^3}\Big ( (2+\ln y^u)(1-y^u)+ 2 y^u \ln y^u \Big
   ).
\end{equation*}
   It is easy to check that the function $h(w)=2(1-w)+(1+w)\ln w$ is concave on $(0,
   1)$. This combined with the observation that $h'(1)=0$ implies
   that $h'(w) \geq 0$ for $0 < w \leq 1$ and hence $h(w) \leq
   h(1)=0$ for $0 < w \leq 1$. Apply this with $w=q^u$ implies that
   $g'(u) \leq 0$. We deduce from this that $f''(x) \geq 0$. Hence
   it remains to show that $f'(1) \geq 0$. In this case we regard
   $f'(1)$ as a function of $s$ with $0<s<1$ and define
\begin{equation*}
  a(s)= \ln \frac {q^{s-1}-1}{(1-s)(1-q)} - \ln q^{s-1}- \frac {\ln
  q^{s-1}}{q^{s-1}-1}+\ln q -\frac {\ln q}{1-q}.
\end{equation*}
   It suffices to show $a(s) \geq 0$ and calculation yields
\begin{equation*}
   a'(s)= \frac {1}{(1-s)(q^{s-1}-1)^2}\Big( (q^{s-1}-1)^2  - (\ln
   q^{s-1})^2q^{s-1} \Big )
\end{equation*}
   If we set $z=q^{s-1}$ so that $z \geq 1$ and define
   $b(z)=(z-1)^2-z(\ln z)^2$, then it is easy to check that $b(z)$
   is convex for $z \geq 1$. As $b'(1)=0$, this implies that $b'(z)
   \geq 0$ for $z \geq 1$ so that $b(z) \geq b(1)=0$ for $z \geq 1$.
   It follows from this that $a'(s) \geq 0$ so that $a(s) \geq
   a(0)=0$ and this completes the proof.
\end{proof}

\begin{lemma}
\label{lem5}
   Let $m>n\geq 1$ be two integers, then for any fixed constant $0<c<1$,
   the function
\begin{equation*}
   a(t;m,n,c)=t^{m-n}+t^n-c(1+t^m)
\end{equation*}
   has exactly one root when $t \geq 1$.
\end{lemma}
\begin{proof}
   We have
\begin{equation*}
   a'(t;m,n,c)=t^{m-1}\Big ((m-n)t^{-n}+nt^{n-m}-cm \Big ).
\end{equation*}
   The function $t \mapsto (m-n)t^{-n}+nt^{n-m}-cm$ is clearly
   decreasing when $t \geq 1$. By considering the cases $t=1$ and $t
   \rightarrow +\infty$ we conclude that $a'(t;m,n,c)$ has exactly one
   root when $t \geq 1$. It follows from this and Cauchy's mean
   value theorem that $a(t;m,n,c)$ has at most two roots when $t\geq 1$.
   This combined with the observation that $a(1;m,n,c)>0$ and $\lim_{t \rightarrow
   +\infty}a(t;m,n,c)<0$ yields the desired conclusion.
\end{proof}

\begin{lemma}
\label{lem6}
   For $t \geq s \geq 0$, we have
\begin{equation*}
   \frac {s}{1-e^{-s}}\cdot\frac {t-s}{1-e^{-(t-s)}} \geq \frac
   t{1-e^{-t}}.
\end{equation*}
\end{lemma}
\begin{proof}
   We write
\begin{equation*}
  f(x)=\frac {x}{1-e^{-x}}
\end{equation*}
   and we observe that
\begin{equation*}
  \Big ( \frac {f'(t)}{f(t)} \Big )'=\frac {e^{-t}\big(t^2-(e^{t/2}-e^{-t/2})^2 \big
  )}{t^2(1-e^{-t})^2} \leq 0,
\end{equation*}
  where the last inequality follows from
\begin{equation*}
  e^{t/2}-e^{-t/2} \geq t, \hspace{0.1in} t \geq 0.
\end{equation*}
  The above inequality can be established by using the Taylor
  expansion of $e^{t/2}-e^{-t/2}$ at $t=0$. We now deduce that
\begin{equation*}
  \frac {f'(t-s)}{f(t-s)}-\frac
  {f'(t)}{f(t)} \geq 0,
\end{equation*}
  for $t \geq s \geq 0$. This implies that the function $t \mapsto \ln f(t-s)-\ln f(t)$ is increasing for $t >s$. Thus we get
\begin{equation*}
   \ln f(s)+\ln f(t-s) - \ln f(t) \geq  \lim_{t \rightarrow s^+}(\ln f(s)+\ln f(t-s) - \ln f(t))
   =0,
\end{equation*}
   which is the desired result.
\end{proof}
\section{Main Results}
\label{sec 4} \setcounter{equation}{0}
\begin{theorem}
\label{thm1}
   Let $a_i$ and $b_i$ $(i=1,\dotsc,n)$ be real numbers
such that $0\le a_1\le\dotsb\le a_n$, $0\le b_1\le\dotsb\le b_n$,
and $\sum_{i=1}^ka_i\le\sum_{i=1}^kb_i$ for $k=1,\dotsc,n$. If
$f''(x)$ is completely monotonic on $(0, +\infty)$, then
\begin{equation*}
   \exp \Biggl (\sum_{i=1}^n \Bigl ( f ( x+a_i )-f(x+b_i)  \Bigr ) \Biggr )
\end{equation*}
   is logarithmically completely monotonic on $(0, +\infty)$.
\end{theorem}
\begin{proof}
   It suffices to show that
\begin{equation*}
   -\sum_{i=1}^n \Big ( f' \left (x+a_i \right )-f' \left (x+b_i \right ) \Big )
\end{equation*}
    is completely monotonic on $(0, +\infty)$ or for $k \geq 1$,
\begin{equation*}
   (-1)^k\sum_{i=1}^nf^{(k)}(x+a_i) \geq (-1)^k\sum_{i=1}^nf^{(k)}(x+b_i).
\end{equation*}
    By Lemma \ref{lem2}, it suffices to show that $(-1)^kf^{(k)}(x)$
    is decreasing and convex on $(0, +\infty)$ or equivalently,
    $(-1)^kf^{(k+1)}(x) \leq 0$ and $(-1)^kf^{(k+2)}(x) \geq 0$ for
    $k \geq 1$. The last two inequalities hold since we assume that
    $f''(x)$ is completely monotonic on $(0, +\infty)$. This completes
    the proof.
\end{proof}

    As a direct consequence of Theorem \ref{thm1}, we now generalize
    a result of Alzer \cite[Theorem 10]{alz1}:
\begin{cor}
\label{cor1}
   Let $a_i$ and $b_i$ $(i=1,\dotsc,n)$ be real numbers
such that $0\le a_1\le\dotsb\le a_n$, $0\le b_1\le\dotsb\le b_n$,
and $\sum_{i=1}^ka_i\le\sum_{i=1}^kb_i$ for $k=1,\dotsc,n$. Then,
\begin{equation*}
 x \mapsto \prod_{i=1}^n\frac{\Gamma_q(x+a_i)}{\Gamma_q(x+b_i)}
\end{equation*}
 is logarithmically completely monotonic on $(0, +\infty)$.
\end{cor}
\begin{proof}
   Apply Theorem \ref{thm1} to $f(x) = \log \Gamma_q(x)$ and note
   that $f''(x)=\psi'_q(x)$ is completely monotonic on $(0, +\infty)$
   and this completes the proof.
\end{proof}

\begin{theorem}
\label{thm2}
   Let $f''(x)$ be completely monotonic on $(0, +\infty)$, then for $0
   \leq s \leq 1$, the functions
\begin{eqnarray*}
   x & \mapsto & \exp \Biggl (-  \Bigl ( f(x+1)-f(x+s)-(1-s)f'(x+\frac {1+s}2) \Bigr )
   \Biggr
   ), \\
   x & \mapsto & \exp \Biggl ( f(x+1)-f(x+s)-\frac {(1-s)}{2}\Bigl (f'(x+1)+f'(x+s) \Bigr )
   \Biggr
   )
\end{eqnarray*}
   are logarithmically completely monotonic on $(0, +\infty)$.
\end{theorem}
\begin{proof}
   We may assume $0 \leq s<1$. We will prove the first assertion and the second one can
   be shown similarly. It suffices to show that
\begin{equation*}
   f'(x+1)-f'(x+s)-(1-s)f''(x+\frac {1+s}2)
\end{equation*}
    is completely monotonic on $(0, +\infty)$ or for $k \geq 1$,
\begin{equation*}
   \frac 1{1-s}\int^{x+1}_{x+s}(-1)^{k+1}f^{(k+1)}(t)dt \geq (-1)^{k+1}f^{(k+1)}(x+\frac {1+s}2).
\end{equation*}
   The last inequality holds by Lemma \ref{lem3} and our assumption
   that $f''(x)$ is completely monotonic on $(0, +\infty)$. This
   completes the proof.
\end{proof}

\begin{cor}
\label{cor2}
   For $0 \leq s \leq 1$, the functions
\begin{eqnarray*}
  x & \mapsto & \frac {\Gamma_q(x+s)}{\Gamma_q(x+1)}\exp \Biggl ( (1-s)\psi_q(x+\frac {1+s}{2})
  \Biggr
   ),  \\
  x & \mapsto & \frac {\Gamma_q(x+1)}{\Gamma_q(x+s)}\exp \Biggl ( -\frac {(1-s)}{2}\Big (\psi_q(x+1)+\psi_q(x+s)
  \Big ) \Biggr
   )
\end{eqnarray*}
 are logarithmically completely monotonic on $(0, +\infty)$.
\end{cor}
\begin{proof}
   Apply Theorem \ref{thm2} to $f(x) = \log \Gamma_q(x)$ and note
   that $f''(x)=\psi'_q(x)$ is completely monotonic on $(0, +\infty)$
   and this completes the proof.
\end{proof}

   By applying Lemma \ref{lem3} to $f(x)=-\psi_q(x)$, we obtain
\begin{theorem}
\label{thm3}
   For positive $x$ and $0 \leq s \leq 1$,
\begin{equation*}
   \exp \Biggl ( \frac {(1-s)}{2} \Big(\psi_q(x+1)+\psi_q(x+s) \Big )
   \Biggr
   ) \leq \frac {\Gamma_q(x+1)}{\Gamma_q(x+s)} \leq \exp \Biggl ( (1-s)\psi_q \big (x+\frac {1+s}{2} \big )
   \Biggr
   ).
\end{equation*}
\end{theorem}

   The upper bound in Theorem \ref{thm3} is due to Ismail and
   Muldoon \cite{I&M}. Our proof here is similar to that of
   Corollary 3 in \cite{Mer}, which asserts for positive $x$,
\begin{equation}
\label{3.10}
   \frac {(1-s)}{2}\Big (\psi(x+1)+\psi(x+s) \Big )
    < \ln \frac {\Gamma(x+1)}{\Gamma(x+s)} < (1-s)\psi(x+\frac
    {1+s}{2}), \hspace{0.1in} 0 < s < 1.
\end{equation}

   We further note the following integral analogue of Theorem 3.5
   in \cite{alz2}:
\begin{equation*}
   \psi \Big(L_0(b,a) \Big ) \leq \frac 1{b-a}\int^b_a\psi(x)dx \leq
   \psi \Big (L_1(b,a) \Big ), \hspace{0.1in} b>a>0.
\end{equation*}
    It follows from this that for positive $x$
    and $0 \leq s \leq 1$,
\begin{equation*}
   \exp \Biggl ( (1-s)\psi \Big (L_0(x+1,x+s) \Big ) \Biggr
   ) \leq \frac {\Gamma(x+1)}{\Gamma(x+s)} \leq \exp \Biggl ( (1-s)\psi \Big (L_1(x+1, x+s) \Big )
   \Biggr
   ).
\end{equation*}
   We
note that the left-hand side inequality above is contained in
\cite[Lemma 1]{E&P1} and \cite[Theorem 4]{E&G&P}. Now by Lemma
\ref{lem4}, observing that $L_{-1}(x+1,
   x+s)=\sqrt{(x+1)(x+s)}$ and $L_{2}(x+1,
   x+s)=x+(1+s)/2$, we obtain
\begin{equation}
\label{3.1}
   \exp \Biggl ( (1-s)\psi \Big (\sqrt{(x+1)(x+s)} \Big ) \Biggr
   ) \leq \frac {\Gamma(x+1)}{\Gamma(x+s)} \leq \exp \Biggl ( (1-s)\psi \Big (x+(1+s)/2 \Big )
   \Biggr
   ).
\end{equation}
    Note that Alzer \cite[Lemma 2.4]{alz2} has shown that $\psi(e^x)$ is strictly concave on ${\mathbb R}$ and it follows that
    $$\psi(x+1)+\psi(x+s) \leq 2\psi \Big ( \sqrt{(x+1)(x+s)} \Big ),$$
     also note that
    $\psi(x)$ is an increasing function on $(0, +\infty)$ and
    $\sqrt{(x+1)(x+s)} \geq x+s^{1/2}$, we
    see that the inequalities in \eqref{3.1} refine the case $q
    \rightarrow 1$ in Theorem \ref{thm3} and the following result of Kershaw \cite{Ker},
    which states that for positive $x$ and $0 \leq s \leq 1$,
\begin{equation*}
   \exp \Big ( (1-s)\psi(x+s^{1/2}) \Big
   ) \leq \frac {\Gamma(x+1)}{\Gamma(x+s)} \leq \exp \Big ( (1-s)\psi(x+\frac {1+s}{2}) \Big
   ).
\end{equation*}

   We now show the lower bound above and the corresponding one
 in \eqref{3.10} are not comparable
   in general (see \cite[p. 856]{E&P} for a similar discussion). In fact, 
   by Theorem 3.7 of \cite{alz2}, we have
\begin{equation*}
   \psi(1)+\psi(s) < 2\psi(s^{1/2}), ~~0<s<1.
\end{equation*}
    This gives one inequality on letting $x \rightarrow 0$. On the other hand, using the well-known series representation (see, for example, \cite[(1.8)]{I&M}):
\begin{equation*}
   \psi(x)=-\gamma+\sum^{\infty}_{n=0}(\frac 1{n+1}-\frac 1{x+n})
\end{equation*}
    with $\gamma=0.57721\ldots$ denoting Euler's constant, we obtain
    for $x>1$,
\begin{equation*}
   \psi(x+1)+\psi(x+s) - 2\psi(x+s^{1/2})=\sum^{\infty}_{n=0}\frac
   {(1-s^{1/2})^2(x+n-s^{1/2})}{(x+n+1)(x+n+s)(x+n+s^{1/2})} > 0.
\end{equation*}

   Our next result refines the left-hand side inequality of \eqref{1.4}:
\begin{theorem}
\label{thm8}
  Let $0<s<1$ and $0<q<1$. Let $u(q,s)$ be defined by \eqref{1.4'} and let the function $g_q(x; s, 1, u(q, s))$ be defined as in Theorem \ref{thm5}. Then $g_q(x; s, 1, u(q, s))$ is
logarithmically completely monotonic on $(0,
  +\infty)$.
\end{theorem}
\begin{proof}
   Define
\begin{equation*}
   h_q(x) = - \ln g_q(x;s,1, u(q, s))=-\ln \Gamma_q(x+1)+ \ln
   \Gamma_q(x+s)+(1-s)\ln \Big(\frac {1-q^{x+u(q, s)}}{1-q} \Big ).
\end{equation*}
   It suffices to show that $h'_q(x)$ is completely monotonic on $(0,
  +\infty)$. We have
\begin{equation*}
   h'_q(x) = - \psi_q(x+1) +
   \psi_q(x+s)-(1-s)\ln q  \frac {q^{x+u(q, s)}} {1-q^{x+u(q, s)}}.
\end{equation*}
   Using the expression \eqref{1.12},
   we can rewrite $h'_q(x)$ as
\begin{equation*}
   h'_q(x) =  -\ln q \Big ( \sum^{\infty}_{n=1}\frac {q^{nx}(q^{n}-q^{ns})}{1-q^n}
   +(1-s) \frac {q^{x+u(q, s)}}{1-q^{x+u(q, s)}} \Big ),  \hspace{0.1in}
   0<q<1.
\end{equation*}
   Expanding
   $(1-q^{x+u(q, s)})^{-1}$, we may further rewrite $h'_q(x)$ as
\begin{equation*}
  h'_q(x) = -\ln q \sum^{\infty}_{n=1}
  \frac {q^{nx}}{1-q^n}w_{q,n}(s),  \hspace{0.1in}
   0<q<1,
\end{equation*}
  where
\begin{equation*}
   w_{q,n}(c)=q^{n}-q^{ns}+(1-s)q^{nu(q, s)}(1-q^n).
\end{equation*}
   In order for $h'_q(x)$ to be completely monotonic on $(0,
  +\infty)$, it suffices to show $w_{q,n}(s) \geq 0$ for $0<s<1$.
  This is just Lemma \ref{lem10} and this completes the proof.
\end{proof}

   The above theorem implies that $g_q(x; s, 1, u(q, s)) > \lim_{x \rightarrow +\infty}g_q(x; s, 1, u(q,
   s))=1$, where the limit can be easily evaluated using \eqref{01}
   and we recover the left-hand side inequality of \eqref{1.4}.

\section{Applications of Lemma \ref{thm6}}
\label{sec 4.1} \setcounter{equation}{0}
   In this section, we look at a few applications of Lemma \ref{thm6} and
   we first deduce Theorem
   \ref{thm3.0} from it.
   As was noted in \cite{G&I}, we have $F_n(x,
   q)=F_{n-1}(x,q)/F_{n-1}(x+a_n,q)$ if $F_n$ and $F_{n-1}$ are
   defined by the same parameters $a_k$. It follows from Lemma
   \ref{thm6} by induction that we only need to prove the case
   $n=1$, which follows from \cite[Lemma 2.1, 3.1]{G&I}.

   We now deduce the following result from Lemma \ref{thm6} and Theorem \ref{thm5}:
\begin{cor}[{\cite[Theorem 4,5]{B&I}}]
\label{cor4}
  The functions
\begin{equation*}
   x \mapsto (1-\frac 1{2x+1/2})^{-1/2}\frac
   {\Gamma^2(x+1/2)}{\Gamma(x)\Gamma(x+1)}, \hspace{0.1in}
   x \mapsto  (1+\frac 1{2x})^{-1/2}\frac {\Gamma(x)\Gamma(x+1)}{\Gamma^2(x+1/2)}
\end{equation*}
  are logarithmically completely monotonic on $(1/4, +\infty)$ and $(0,
  +\infty)$, respectively.
\end{cor}
\begin{proof}
    We take $q=1$, $a=0, c=-1/4, b=1/2$ and $a=c=0, b=1/2$ in Theorem \ref{thm5} to conclude that the
    functions
\begin{equation*}
  f_1(x) =   (x-1/4)^{-1/2}\frac {\Gamma(x+1/2)}{\Gamma(x)}, \hspace{0.1in}
  f_2(x) =  \frac {\Gamma(x)\sqrt{x}}{\Gamma(x+1/2)}
\end{equation*}
    are logarithmically completely monotonic on $(1/4, +\infty)$ and $(0, +\infty)$, respectively.
    The theorem now follows by applying Lemma \ref{thm6} to the
    functions $f_1(x)/f_1(x+1/2)$ and $f_2(x)/f_2(x+1/2)$.
\end{proof}
   We remark here the above corollary is essentially Theorem 4 and 5 in
\cite{G&I}, except that Theorem 4 of \cite{B&I} is originally stated
as the function
\begin{equation*}
   x \mapsto (1-\frac 1{2x})^{-1/2}\frac
   {\Gamma^2(x+1/2)}{\Gamma(x)\Gamma(x+1)}
\end{equation*}
    is completely monotonic on $(1/2, +\infty)$. One
    can also obtain this result by modifying our approach above.

     As another application of Lemma \ref{thm6}, we note by \cite[Theorem 2.1]{ILM},
     the function
\begin{equation*}
   f(x; \alpha) = x^{\alpha}\Gamma(x)(\frac {e}{x})^x
\end{equation*}
    is logarithmically completely monotonic on $(0, +\infty)$ for $\alpha \leq 1/2$, and
    $1/f(x; \alpha)$ is logarithmically completely monotonic on $(0, +\infty)$ for $\alpha \geq 1$.
    It
    follows from this and Lemma \ref{thm6} that for any $a>0$, the functions
    $f(x;1/2)/f(x+a;1/2)$ and $f(x+a;1)/f(x;1)$ are logarithmically completely monotonic on $(0,
    +\infty)$. Hence by \eqref{3.05},
\begin{equation*}
   \frac {f(x;1/2)}{f(x+a;1/2)} \geq \lim_{x \rightarrow +\infty} \frac
   {f(x;1/2)}{f(x+a;1/2)}=1, \hspace{0.1in} \frac {f(x+a;1)}{f(x;1)} \geq \lim_{x \rightarrow +\infty} \frac
   {f(x+a;1)}{f(x;1)}=1.
\end{equation*}
    We can recast the above as
\begin{equation*}
   \frac {b^{b-1}}{a^{a-1}}e^{a-b} \leq \frac {\Gamma(b)}{\Gamma(a)}
   \leq \frac {b^{b-1/2}}{a^{a-1/2}}e^{a-b}, \hspace{0.1in} b>a>0.
\end{equation*}
    The above inequalities are due to Ke\v cki\'c and Vasi\'c \cite{K&V} for the case $b \geq a >1$.

   We now look at some functions involving the $q$-gamma function. Ismail et al.
   \cite[Theorem 2.2]{ILM} proved that for $q \in (0, 1)$ the function
   $(1-q)^x\Gamma_q(x)$ is logarithmically completely monotonic on $(0,
   +\infty)$. It follows from this and Lemma \ref{thm6} that
   $(1-q)^x\Gamma_q(x)/(1-q)^{x+a}\Gamma_q(x+a)$ is also
   logarithmically completely monotonic on $(0, +\infty)$ for any
   $a>0$. Apply Lemma \ref{thm6} again to the above function, we
   conclude that the function
\begin{equation*}
  g(x)=  \frac {(1-q)^x\Gamma_q(x)(1-q)^{x+2a}\Gamma_q(x+2a)}{((1-q)^{x+a}\Gamma_q(x+a))^2}
\end{equation*}
  is logarithmically completely monotonic on $(0, +\infty)$ for any
  $a>0$. It follows from this by using \eqref{01} that $g(x) \geq \lim_{x \rightarrow
  +\infty}g(x)=1$. Other than the equality sign, this is a result of Alzer and Grinshpan \cite[Theorem
  3.1]{A&G}. One can obtain a similar result for $q>1$.

   In \cite[Lemma 2.7]{A&G}, Alzer and Grinshpan generalized a
   result of Ismail and Muldoon \cite{I&M} that for $q>0$, the
   function $x \mapsto \Gamma_q(x)/\Gamma^{1/\beta}_{q^{1/\beta}}(\beta x)$ is logarithmically completely monotonic on $(0,
   +\infty)$ for $\beta >1$ and the
   function $x \mapsto\Gamma^{1/\beta}_{q^{1/\beta}}(\beta x)/ \Gamma_q(x)$ is logarithmically completely monotonic on $(0,
   +\infty)$ for $0 < \beta <1$. Proceeding similarly to our
   discussions above, we see that the function
\begin{equation*}
   z \mapsto \frac {\Gamma_q(z)\Gamma_q(x+y+z)
   \Gamma^{1/\beta}_{q^{1/\beta}}(\beta(z+x))\Gamma^{1/\beta}_{q^{1/\beta}}(\beta (z+y))}
   {\Gamma_q(z+x)\Gamma_q(z+y)\Gamma^{1/\beta}_{q^{1/\beta}}(\beta z)\Gamma^{1/\beta}_{q^{1/\beta}}(\beta (z+x+y))}
\end{equation*}
   is is logarithmically completely monotonic on $(0,
   +\infty)$ for $\beta >1, x>0, y>0$ and a similar assertion holds for $0 < \beta <1$. We deduce readily from this the
   following
\begin{cor}[{\cite[Theorem 2.8]{A&G}}]
\label{cor5}
  Let $q>0$ and $x, y, z>0$, then for any $\alpha >1$, the
  inequality
\begin{equation*}
  \Big ( \frac {\Gamma_{q^{\alpha}}(z+x)\Gamma_{q^{\alpha}}(z+y)}
  {\Gamma_{q^{\alpha}}(x+y+z)\Gamma_{q^{\alpha}}(z)} \Big )^{\alpha} \leq \frac {\Gamma_q(\alpha(z+x))\Gamma_q(\alpha(z+y))}{\Gamma_q(\alpha z)\Gamma_q(\alpha(x+y+z))}
\end{equation*}
  holds. If $\alpha \in (0, 1)$, then the reversed inequality holds.
\end{cor}

  For $\alpha \geq 0$, let
\begin{equation}
\label{4.1}
  f_{\alpha}(x)=-\ln \Gamma(x)+(x-\frac 1{2})\ln x -x +\frac {1}{12}\phi'(x+\alpha).
\end{equation} 
  Alzer \cite[Theorem 1]{alz0} showed that for fixed $0<s<1$, the function $\exp(-f_{\alpha}(x+s)/f_{\alpha}(x+1))$ is strictly completely monotonic on $(0, +\infty)$ if and only if $\alpha \geq 1/2$ and the function $\exp(f_{\alpha}(x+s)/f_{\alpha}(x+1))$ is strictly completely monotonic on $(0, +\infty)$ if and only if $\alpha =0$. In view of Lemma \ref{thm6}, it suffices to show that $f'_{\alpha}(x)$ is strictly completely monotonic on $(0, +\infty)$ if and only if $\alpha \geq 1/2$ and that $-f'_{\alpha}(x)$ is strictly completely monotonic on $(0, +\infty)$ if and only if $\alpha =0$. We now prove the ``if" part of the above assertions in the following
\begin{prop}
\label{prop4.1}
  For $\alpha \geq 0$, let $f_{\alpha}(x)$ be defined as in \eqref{4.1}. Then $f'_{\alpha}(x)$ is strictly completely monotonic on $(0, +\infty)$ if $\alpha \geq 1/2$ and that $-f'_{\alpha}(x)$ is strictly completely monotonic on $(0, +\infty)$ if $\alpha =0$.
\end{prop}
\begin{proof}
  As $\phi'(x)$ is completely monotonic on $(0, +\infty)$, we may just focus on the case $\alpha=1/2$ or $0$. We have
\begin{equation*}
  f'_{\alpha}(x)=-\psi(x)+\ln x -\frac 1{2x}+\frac 1{12}\phi''(x+\alpha),  \hspace{0.1in} f''_{\alpha}(x)=-\psi'(x)+\frac 1{x} +\frac 1{2x^2}+\frac 1{12}\phi^{(3)}(x+\alpha).
\end{equation*}
  Using the asymptotic expressions \eqref{3.03} and \eqref{1.2}, we see that $\lim_{x \rightarrow +\infty}f^{(n)}_{\alpha}(x)=0$ for any integer $n \geq 0$. It is then easy the see that the assertions of the proposition will follows if we can show that $f''_{1/2}(x+1)-f''_{1/2}(x)$ is strictly completely monotonic on $(0, +\infty)$ and $f''_{0}(x)-f''_{0}(x+1)$  is completely monotonic on $(0, +\infty)$. Using \eqref{3.4}, it is easy to see that
\begin{equation*}
  f''_{1/2}(x+1)-f''_{1/2}(x)=\frac {x^2+x+1/8}{4x^2(1+x)^2(x+1/2)^4}=\frac {1}{4(1+x)^2(x+1/2)^4}(1+\frac {1}{x}+\frac {1}{8x^2}).
\end{equation*}
   It is easy to see from this and the ring structure of completely
  monotonic functions on $(0, +\infty)$ that $f''_{1/2}(x+1)-f''_{1/2}(x)$ is strictly completely monotonic on $(0, +\infty)$. Similarly, one shows that $f''_{0}(x)-f''_{0}(x+1)$  is completely monotonic on $(0, +\infty)$ and this completes the proof.
\end{proof}

   We note here that recently, Alzer and Batir \cite{AB} showed that the following function ($x > 0, c \geq 0$)
\begin{equation*}
   G_{c}(x)=\ln \Gamma(x)-x\ln x+x-\frac 1{2} \ln (2\pi)+\frac 1{2}\psi(x+c)
\end{equation*}
   is completely monotonic if and only if $c \geq 1/3$ and $-G_c(x)$ is completely monotonic if and only if $c = 0$. We point out here that using the idea in the proof of Proposition \ref{prop4.1}, one can give another proof of the above result of Alzer and Batir and we shall leave the details to the reader.
\section{Inequalities Involving Polygamma Functions }
\label{sec 4.2} \setcounter{equation}{0}
  For integers $p \geq m \geq n \geq q \geq 0$ and any real number $c$, we
   define
\begin{equation*}
   F_{p, m,n,q}(x; c)= (-1)^{m+n}\psi^{(m)}(x)\psi^{(n)}(x) - c(-1)^{p+q}\psi^{(p)}(x)\psi^{(q)}(x).
\end{equation*}
   Here we set $\psi^{(0)}(x)=-1$ for convenience.

   For $n \geq 2$, a result of Alzer and Wells \cite[Theorem
2.1]{A&W} asserts that
   the function $F_{n+1, n, n, n-1}(x; c)$ is completely monotonic on $(0, +\infty)$ if and only if $c \leq
   (n-1)/n$ and $-F_{n+1, n, n, n-1}(x; c)$ is completely monotonic on $(0, +\infty)$ if and only if $c \geq
   n/(n+1)$.

   We denote
\begin{equation*}
 c_{p,m,n,q} =
     \frac {(m-1)!(n-1)!}{(p-1)!(q-1)!}, \hspace{0.1in}  q \geq 1; \hspace{0.1in}
   c_{p,m,n,0} =  \frac {(m-1)!(n-1)!}{(p-1)!},
\end{equation*}
   and
\begin{equation*}
   d_{p,m,n,q} =
     \frac {m!n!}{p!q!}.
\end{equation*}
   and note that $0<c_{p,m,n,q}, d_{p,m,n,q}<1$ when $p+q=m+n, p>m$. We now extend the result of Alzer and Wells to the following
\begin{theorem}
\label{thm4}
   Let $p>m \geq n >q \geq 0$ be integers satisfying $m+n=p+q$. The
   function
   $F_{p, m,n,q}(x; c_{p,m,n,q})$ is completely monotonic on $(0, +\infty)$. The function $-F_{p, m,n,q}(x; d_{p, m,n, q})$
   is also completely monotonic on $(0, +\infty)$ when $q>0$, .
\end{theorem}
\begin{proof}
  We first prove the assertion for $F_{p, m,n,q}(x; c_{p,m,n,q})$ with
  $q \geq 1$ and the proof here uses the method in \cite{A&W}.
  Using the integral
  representation \eqref{1.10} for $(-1)^{n+1}\psi^{(n)}(x)$ and using $*$ for the
  Laplace convolution, we get
\begin{equation*}
   F_{p, m,n,q}(x; c_{p,m,n,q})=\int^{\infty}_0e^{-xt}g(t)dt,
\end{equation*}
   where
\begin{eqnarray*}
   g(t) &=& \frac {t^m}{1-e^{-t}} * \frac {t^n}{1-e^{-t}}-c_{p,m,n,q}\frac
{t^p}{1-e^{-t}} * \frac {t^q}{1-e^{-t}} \\
  &=& \int^t_{0}\Big ( (t-s)^ms^n-c_{p,m,n,q}(t-s)^ps^q \Big )h(t-s)h(s)ds,
\end{eqnarray*}
  with
\begin{equation*}
   h(s)=\frac 1{1-e^{-s}}.
\end{equation*}
   It suffices to show $g(t) \geq 0$ and by a change of variable
   $s \rightarrow ts$ we can recast it as
\begin{equation*}
   g(t)=t^{m+n+1}\int^1_{0}\Big ( (1-s)^ms^n-c_{p,m,n,q}(1-s)^ps^q \Big )
   h \Big(t(1-s) \Big)h(ts)ds.
\end{equation*}
   We now break the above integral into two integrals, one from $0$
   to $1/2$ and the other from $1/2$ to $1$. We
   make a further change of variable $s \rightarrow (1-s)/2$ for the first one and
   $s \rightarrow (1+s)/2$ for the second one. We now combine them
   to get
\begin{equation*}
  g(t)=\Big(\frac {t}2 \Big )^{m+n+1}\int^1_{0}a(\frac {1+s}{1-s}; p-q, n-q, c_{p,m,n,q})
  (1-s^2)^q(1-s)^{p-q}
   h \Big ( \frac {t(1-s)}2 \Big )h\Big(\frac {t(1+s)}2\Big )ds,
\end{equation*}
   where the function $a(t;m,n,c)$ is as defined in Lemma
   \ref{lem5}. Note that $(1+s)/(1-s) \geq 1$ for $0 \leq s <1$,
   hence by Lemma \ref{lem5}, there is a unique number $0 < s_0< 1$ such
   that
\begin{equation*}
   a \Big ( \frac {1+s_0}{1-s_0}; p-q, n-q, c_{p,m,n,q} \Big)=0.
\end{equation*}
   We further note it is shown in the proof of \cite[Lemma 2.2]{A&W}
   that the function
\begin{equation*}
   s \mapsto (1-s^2)h \Big ( \frac {t(1-s)}2 \Big )h\Big(\frac {t(1+s)}2\Big )
\end{equation*}
   is a decreasing function on $(0,1)$ so that for $0 \leq s \leq
   1$,
\begin{eqnarray*}
   && a(\frac {1+s}{1-s}; p-q, n-q, c_{p,m,n,q})
  (1-s^2)^q(1-s)^{p-q}
   h \Big ( \frac {t(1-s)}2 \Big )h\Big(\frac {t(1+s)}2\Big ) \\
   &\geq & a(\frac {1+s}{1-s}; p-q, n-q, c_{p,m,n,q})
  (1-s^2)^{q-1}(1-s)^{p-q}
  (1-s^2_0) h \Big ( \frac {t(1-s_0)}2 \Big )h\Big(\frac {t(1+s_0)}2\Big
  ).
\end{eqnarray*}
   Hence
\begin{eqnarray*}
   g(t) &\geq& \Big(\frac {t}2 \Big )^{m+n+1}(1-s^2_0) h \Big ( \frac {t(1-s_0)}2 \Big )h\Big(\frac {t(1+s_0)}2\Big
  ) \\
  && \cdot \int^1_{0}a(\frac {1+s}{1-s}; p-q, n-q, c_{p,m,n,q})
  (1-s^2)^{q-1}(1-s)^{p-q}ds.
\end{eqnarray*}
  Note that the integral above is (by reversing the process above on
  changing variables)
\begin{equation*}
   2^{m+n-1}\int^1_{0}\Big ( (1-s)^{m-1}s^{n-1}-c_{p,m,n,q}(1-s)^{p-1}s^{q-1} \Big
   )ds=0,
\end{equation*}
   where the last step follows from the well-known beta function
   identity
\begin{equation*}
   B(x,y)=\int^1_0t^{x-1}(1-t)^{y-1}dt=\frac
   {\Gamma(x)\Gamma(y)}{\Gamma(x+y)}, \hspace {0.1in} x,y>0,
\end{equation*}
   and the well-known fact $\Gamma(n) = (n-1)!$ for $n \geq 1$.

   Now we prove the assertion for $F_{p, m,n,0}(x; c_{p,m,n,0})$. In this case $p=m+n$
   and we note that
\begin{equation*}
   c_{m+n,m,n, 0} = B(m,n)=\int^1_0s^{m-1}(1-s)^{n-1}ds,
\end{equation*}
   and we use this to write
\begin{equation*}
   c_{m+n,m,n,0}\frac {t^{m+n}}{1-e^{-t}}=\int^t_0\frac
   {s^{m-1}(t-s)^{n-1}t}{1-e^{-t}}ds.
\end{equation*}
   It follows that
\begin{eqnarray*}
   && F_{p, m,n,0}(x; c_{p,m,n,0}) \\
   &=& \int^{\infty}_0e^{-xt}
   \Big ( \frac {t^m}{1-e^{-t}} * \frac {t^n}{1-e^{-t}} - c_{m+n,m,n,0}\frac {t^{m+n}}{1-e^{-t}} \Big
   )dt \\
   &=&  \int^{\infty}_0e^{-xt} \left ( \int^t_0s^{m-1}(t-s)^{n-1}
   \left ( \frac {s}{1-e^{-s}}\cdot\frac {t-s}{1-e^{-(t-s)}}-
   \frac t{1-e^{-t}}
   \right )ds \right ) dt \\
   & \geq & 0,
\end{eqnarray*}
   where the last inequality follows from Lemma \ref{lem6}.

   It remains to show the assertion for $-F_{p, m,n,q}(x; d_{p, m,n,
   q})$. In this case we use the series representation in
   \eqref{1.10} for $(-1)^{n+1}\psi^{(n)}(x)$ to get
\begin{eqnarray*}
   && -F_{p, m,n,q}(x; d_{p, m,n, q}) \\
   &=& m!n!\left ( \left(\sum^{\infty}_{i=0}\frac 1{(x+i)^{p+1}}
   \right)\left(\sum^{\infty}_{j=0}\frac 1{(x+j)^{q+1}}\right)-\left(\sum^{\infty}_{i=0}\frac 1{(x+i)^{m+1}}
   \right)\left(\sum^{\infty}_{j=0}\frac 1{(x+j)^{n+1}}\right)\right
   ).
\end{eqnarray*}
   We note the following Binet-Cauchy identity:
\begin{equation*}
   (\sum^n_{i=0}a_ic_i)(\sum^n_{i=0}b_id_i)-(\sum^n_{i=0}a_id_i)(\sum^n_{i=0}b_ic_i)=\sum_{0 \leq
   i<j\leq n}(a_ib_j-a_jb_i)(c_id_j-c_jd_i).
\end{equation*}
   We now apply the above identity with
\begin{equation*}
   a_i=\frac 1{(x+i)^{m+1}}, b_i=\frac 1{(x+i)^{q+1}}, c_i=\frac
   1{(x+i)^{p-m}}, d_i=1
\end{equation*}
   to get
\begin{eqnarray*}
   && -F_{p, m,n,q}(x; d_{p, m,n, q}) \\
   &=& m!n!\sum_{0 \leq i<j}
   \left ( \frac 1{(x+i)^{m+1}(x+j)^{q+1}}-\frac 1{(x+j)^{m+1}(x+i)^{q+1}}
   \right)\left ( \frac 1{(x+i)^{p-m}}-\frac 1{(x+j)^{p-m}} \right)
\end{eqnarray*}
   We note that the second factor on the right-hand side above is completely monotonic
   on $(0, +\infty)$ and also that
\begin{eqnarray*}
   &&\frac 1{(x+i)^{m+1}(x+j)^{q+1}}-\frac 1{(x+j)^{m+1}(x+i)^{q+1}} \\
   &=& \frac 1{(x+i)^{q+1}}\frac 1{(x+j)^{q+1}}
   \Big ( \frac 1{(x+i)^{m-q}}-\frac 1{(x+j)^{m-q}}\Big ).
\end{eqnarray*}
  Certainly each factor on the right-hand side above is completely monotonic on $(0,
  +\infty)$ and it follows from the ring structure of completely
  monotonic functions on $(0, +\infty)$ that the left-hand side
  above is also completely monotonic on $(0, +\infty)$. Hence by the ring structure of completely
  monotonic functions on $(0, +\infty)$ again we deduce that $-F_{p, m,n,q}(x; d_{p, m,n,
  q})$ is completely monotonic on $(0, +\infty)$.
\end{proof}
   We point out here corresponding to $m=n=1, p=2, q=0$ in Theorem
   \ref{thm4}, it was shown in the
   proof of \cite[(4.39)]{alz2} and in \cite[Lemma 1.1]{Ba} the following special case (in fact with
   strict inequality)
\begin{equation}
\label{4.2}
   \Big ( \psi'(x) \Big )^2+\psi''(x) \geq 0, \hspace{0.1in} x>0.
\end{equation}
   We now give another proof of the above inequality by establishing the following
\begin{prop}
\label{prop5.1}
 Let $0<c<1$ be fixed. Then for any $x>0$,
\begin{equation}
\label{5.2}
  \frac 1{c}\Big(\psi(x+c)-\psi(x) \Big )^2 > \psi'(x)-\psi'(x+c) > \Big(\psi(x+c)-\psi(x) \Big )^2.
\end{equation}
  The above inequalities reverse when $c>1$.
\end{prop}
\begin{proof}
  We first prove the left-hand side inequality of \eqref{5.2}. For this, we define
\begin{equation*}
  f(x)= \psi'(x+c)-\psi'(x)+\frac 1{c}\Big(\psi(x+c)-\psi(x) \Big )^2.
\end{equation*}
  Applying the relations \eqref{3.4}, we obtain
\begin{equation*}
  f(x+1)-f(x)= -c\Big (\frac 1{x}-\frac 1{x+c} \Big ) \Big( 2\psi(x+c)-2\psi(x)+\frac {1-c}{x+c}-\frac {1+c}{x} \Big ).
\end{equation*}
  We now apply \eqref{1.10} to rewrite the last expression above as
\begin{equation*}
  2\psi(x+c)-2\psi(x)+\frac {1-c}{x+c}-\frac {1+c}{x}=\int^{\infty}_{0}\frac {e^{-xt}}{1-e^{-t}}\Big ( 2(1-e^{-ct})+((1-c)e^{-ct}-1-c)(1-e^{-t}) \Big ) dt.
\end{equation*}
  We now define for $t>0$,
\begin{equation*}
   g(t)= 2(1-e^{-ct})+((1-c)e^{-ct}-1-c)(1-e^{-t}).
\end{equation*}
  Then we have
\begin{equation*}
  g'(t)=(1+c)e^{-t}(ce^{(1-c)t}+(1-c)e^{-ct}-1) \geq 0,
\end{equation*} 
  when $0<c<1$, where the last inequality follows from applying the arithmetic-geometric inequality to $ce^{(1-c)t}+(1-c)e^{-ct}$ with equality holding if and only if $c=0$ or $1$ or $t=0$. Similarly $g'(t) \leq 0$ when $c>1$. From now on we will assume $0<c<1$ and the case $c>1$ can be discussed similarly. We conclude that $g(t) > g(0)=0$ for $t>0$ when $0<c<1$ and it follows that $f(x+1)-f(x)<0$ in this case. As $\lim_{x \rightarrow +\infty}f(x)=0$, we conclude that $f(x)>0$ which completes the proof of the left-hand side inequality of \eqref{5.2}. The proof of the right-hand side inequality of \eqref{5.2} follows similarly except it is easier and we shall leave it to the reader.
\end{proof}

  We note here the left-hand side inequality of \eqref{5.2} was established in the proof of \cite[Theorem
   1.1]{C} and we rewrite it as
\begin{equation*}
   \frac {\psi'(x+c)-\psi'(x)}{\psi(x+c)-\psi(x)}+\frac {\psi(x+c)-\psi(x)}{c} > 0,
\end{equation*}
   for $x>0, 0 < c< 1$. From this we get back \eqref{4.2} via taking the limit $c
   \rightarrow 0$. We note also that the right-hand side inequality was proven in \cite[Lemma 7]{alz10} and our proof here is simpler for both inequalities of \eqref{5.2}.

    We now prove a $q$-analogue to Proposition \ref{prop5.1}:
\begin{theorem}
\label{thm5.2}
 Let $0<q<1$ and $0<c<1$ be fixed. Then for any $x>0$,
\begin{equation}
\label{5.3}
  \frac {1-q}{1-q^c}\Big(\psi_q(x+c)-\psi_q(x) \Big )^2 > q^x(\psi'_q(x)-\psi'_q(x+c))> \Big(\psi_q(x+c)-\psi_q(x) \Big )^2.
\end{equation}
  The above inequalities reverse when $c>1$.
\end{theorem}
\begin{proof}
  We first prove the left-hand side inequality of \eqref{5.3}. For this, we define
\begin{equation*}
  f(x)= q^x(\psi'_q(x+c)-\psi'_q(x))+\frac {1-q}{1-q^c}\Big(\psi_q(x+c)-\psi_q(x) \Big )^2.
\end{equation*}
  Applying \eqref{1.12}, we obtain
\begin{eqnarray*}
  f(x) &=& (\ln q )^2  \sum^{\infty}_{n=1}\frac {nq^{(n+1)x}(q^c-1)}{1-q^n}+\frac {1-q}{1-q^c}\Big( \ln q \sum^{\infty}_{n=1}\frac {q^{nx}(q^{nc}-1)}{1-q^n }\Big )^2 \\
  &=& \frac {(1-q)(\ln q )^2}{1-q^c}  \sum^{\infty}_{n=3} q^{nx}g_n(c;q), 
\end{eqnarray*}
  where
\begin{equation*}
  g_{n}(c;q)= \sum^{n-1}_{k=1}\Big ( \frac {(1-q^{kc})(1-q^{(n-k)c})}{(1-q^{k})(1-q^{n-k})}-\frac {(1-q^{c})(1-q^{(n-1)c})}{(1-q)(1-q^{n-1})} \Big )
\end{equation*}
   It suffices to show that $g_c(x;q) > 0$ for $0<c<1$ and $g_cx;q) <0$ for $c>1$ when $n \geq 3$. For this, we let $y=q^c$ so that $0<y<1$ and it suffices to show the function $(1-y^k)(1-y^{n-k})/((1-y)(1-y^{n-1}))$ is increasing for $0<y<1$ and $1 \leq k \leq n-1$. On taking the logarithmic derivative of the above function, we see that it suffices to show that $h(k;y)+h(n-k;y) \leq h(1;y)+h(n-1;y)$, where
\begin{equation*}
  h(z;y)=\frac {z}{1-y^z}.
\end{equation*}
  We now regard $h(z;y)$ as a function of $z$ and note that
\begin{equation*}
  h''(z;y)=\frac {2(\ln y)y^z}{(1-y^z)^3}u(y^z), \hspace{0.1in} u(t)=2-2t+\ln t+t\ln t.
\end{equation*}
  It is easy to see that $u'(t)>0$ for $0<t<1$ so that $u(t) < u(1)=0$ for $0<t<1$. It follows that $h''(z;y)>0$. We then deduce from this and Lemma \ref{lem2} that $h(k;y)+h(n-k;y) \leq h(1;y)+h(n-1;y)$ holds and the completes the proof for the left-hand side inequality of \eqref{5.3}.
  For the right-hand side inequality of \eqref{5.3}, one proceeds similarly to the above argument to see that it suffices to show the function $(1-y^k)(1-y^{n-k})/(1-y^{n-1})$ is decreasing for $0<y<1$ and $1 \leq k \leq n-1$. This follows from the observation that both functions $(1-y^k)/(1-y^{n-1})$ and $1-y^{n-k}$ are decreasing and this completes the proof.
\end{proof}
 
   Similar to our discussion above, we have
\begin{cor}
\label{cor5.1}
   For $0<q<1$ and $x>0$,
\begin{equation*}
   \Big ( \psi'_q(x) \Big )^2+\frac {(\ln \frac 1{q}) q^x}{1-q}\psi''_q(x) \geq 0.
\end{equation*}
\end{cor}


   A result of Alzer \cite[Lemma 2]{alz1.6} (see also \cite[Lemma 2.1]{alz2}) asserts that the function
\begin{equation}
\label{4.3}
   x \mapsto -\frac {x\psi^{(n+1)}(x)}{\psi^{(n)}(x)}
\end{equation}
   is strictly decreasing from $[0, +\infty)$ onto
   $(n,n+1]$ for any integer $n \geq 1$. This result implies
   the second assertion in Lemma \ref{thm11} and that the function
   $x \mapsto x^{n+1}(-1)^{n+1}\psi^{(n)}(x)$ is strictly increasing
   for $x>0$.
   In a similar fashion, we
   have
\begin{prop}
\label{cor4.5}
  Let $a \geq 1/2$ be fixed, the for any integer $n \geq 1$, the function
\begin{equation*}
  r_{n,a}(x)=\frac {x\psi^{(n+1)}(x+a)}{\psi^{(n)}(x+a)}
\end{equation*}
   is a decreasing function for $x \geq 0$.
\end{prop}
\begin{proof}
  We have
\begin{eqnarray*}
  &&\Big(\psi^{(n)}(x+a) \Big )^2r'_{n,a}(x) \\
 &=& \psi^{(n+1)}(x+a)\Big(\psi^{(n)}(x+a)+\frac {x}{n}\psi^{(n+1)}(x+a) \Big
 ) \\
 && + x\Big(\psi^{(n)}(x+a)\psi^{(n+2)}(x+a)-(1+\frac {1}{n})\big ( \psi^{(n+1)}(x+a) \big )^2 \Big
 ).
\end{eqnarray*}
  It follows from Lemma \ref{thm11} that
\begin{equation*}
  \psi^{(n+1)}(x+a)\Big(\psi^{(n)}(x+a)+\frac {x}{n}\psi^{(n+1)}(x+a) \Big
 ) \leq 0,
\end{equation*}
  and from Theorem \ref{thm4} that
\begin{equation*}
   \psi^{(n)}(x+a)\psi^{(n+2)}(x+a) \leq \frac {n+1}{n}\big ( \psi^{(n+1)}(x+a) \big )^2.
\end{equation*}
  Hence we conclude that $r'_{n,a}(x) \leq 0$ and this completes the
  proof.
\end{proof}
  We note here the above result implies that
  $r_{n,a}(x) \geq \lim_{x \rightarrow +\infty}r_{n,a}(x)=-n$, which gives back the first assertion in
  Lemma \ref{thm11}.

  One may ask whether it is true or not that the function defined in
  \eqref{4.3} is completely monotonic on $(0, +\infty)$. We note
  here this is not true in general:
  if the function defined in
  \eqref{4.3} is completely monotonic on $(0, +\infty)$ for any
  integer $n \geq 1$, it will follow by induction and the fact that
  $x\psi'(x)$ is completely monotonic on $(0, +\infty)$ (see the paragraph below Lemma \ref{thm11})
   and the ring structure of completely monotonic functions that $(-1)^{n+1}x^n\psi^{(n)}(x)$
   is completely monotonic on $(0,
   +\infty)$. This will further imply a conjecture of Clark and
   Ismail \cite{C&M}, which asserts that the $n$-th derivative of
   $(-1)^{n+1}x^n\psi^{(n)}(x)$ is completely monotonic on $(0,
   +\infty)$ and this has been recently disproved by Alzer et al.
   \cite{ABK}.

  Using the integral representation for $(-1)^{n+1}\psi^{(n)}(x)$ in \eqref{1.10},
  we obtain via integration by parts that
\begin{equation*}
   (-1)^{n+1}x\psi^{(n)}(x)=\int^{\infty}_0\frac
   {t^n}{1-e^{-t}}d(-e^{-xt})=\frac {-e^{-xt}t^n}{1-e^{-t}}\Big
   |^{+\infty}_0+\int^{\infty}_0e^{-xt}\Big ( \frac {d}{dt}\frac
   {t^n}{1-e^{-t}} \Big )dt.
\end{equation*}
   Note that
\begin{equation*}
   \frac {d^m}{dt^m}\frac {t^n}{1-e^{-t}}\Big
   |_{t=0}=0
\end{equation*}
   for $0 \leq m \leq n-2$. Hence repeating the above process, we
   get for $n \geq 1$,
\begin{equation*}
  f_{0,n}(x)=(-1)^{n+1}x^n\psi^{(n)}(x)=c_n+\int^{\infty}_0e^{-xt}\Big ( \frac {d^n}{dt^n}\frac
   {t^n}{1-e^{-t}} \Big )dt,
\end{equation*}
   for some constant $c_n$, where $f_{0,n}(x)$ is defined as in
   Lemma \ref{thm11}. It is shown by Clark and Ismail in \cite{C&M} that
\begin{equation*}
   \frac {d^n}{dt^n}\frac
   {t^n}{1-e^{-t}}
\end{equation*}
   is positive on $(0, +\infty)$ for $2 \leq n \leq 16$ and one checks easily that this is
  also the case for $n=1$. This combined with the observation that
  $f_{0,n}(x)>0$ by \eqref{1.10} and $f'_{0,n}(x)<0$ by Lemma
  \ref{thm11} implies that $f_{0,n}(x)$ is completely monotonic on $(0,
   +\infty)$ for $1 \leq n \leq 16$.

  We further note that similar to Lemma \ref{thm11}, Alzer has shown
  that \cite[Lemma 1]{alz1.6} that the function $xf_{0, n}(x)$ is
  increasing on $(0, +\infty)$ for all $n \geq 1$, where $f_{0,n}(x)$ is defined as in Lemma
  \ref{thm11}. Moreover, he has shown \cite[Lemma 3]{alz1.55} that
  $(xf_{0,1}(x))''$ is strictly completely monotonic on $(0, +\infty)$.
  We note that by our discussion above,
\begin{equation*}
   xf_{0,1}(x)=c_1x+d_1+\int^{\infty}_0e^{-xt}\Big ( \frac {d^2}{dt^2}\frac
   {t}{1-e^{-t}} \Big )dt,
\end{equation*}
   with $c_1, d_1$ constants and $c_1$ as previously defined. One checks
   easily that
   $(t/(1-e^{-t}))'' >0$
   for $t>0$ and this gives another proof of Alzer's result
   mentioned above.

  We now prove a $q$-analogue to Lemma \ref{thm11}:
\begin{theorem}
  For fixed $n \geq 1, 0<q<1$, the function
$f_{n}(x;q)=(1-q^x)^n(-1)^{n+1}\psi^{(n)}_q(x)$ is decreasing on $(0,
+\infty)$.
\end{theorem}
\begin{proof}
  Differentiation gives
\begin{equation*}
  f'_n(x;q)=-(1-q^x)^n\Big ( \frac {n q^x \ln q}{1-q^x}(-1)^{n+1}\psi^{(n)}_q(x)+
  (-1)^{n+2}\psi^{(n+1)}_q(x) \Big ).
\end{equation*}
  It follows from  \eqref{1.12} that for $n \geq 1, 0 <q <1$,
\begin{equation*}
  (-1)^{n+1}\psi^{(n)}_q(x) = (-\ln q )^{n+1} \sum^{\infty}_{k=1}\frac
  {k^nq^{kx}}{1-q^k}.
\end{equation*}
  It follows from this that
\begin{equation*}
  f'_n(x;q)= (1-q^x)^n (-\ln q )^{n+2} \sum^{\infty}_{k=1} \frac
  {q^{kx}}{1-q^k}(n\sum^{k-1}_{m=1}\frac {m^n(1-q^k)}{1-q^m}-k^{n+1}),
\end{equation*}
  where we define the empty sum to be $0$. As the function $q
  \mapsto (1-q^k)/(1-q^m)$ is a decreasing function for $0<q<1$
  with $1 \leq k \leq m$. We conclude that for $0<q<1$,
\begin{equation*}
  n\sum^{k-1}_{m=1}\frac {m^n(1-q^k)}{1-q^m}-k^{n+1} \leq \lim_{q \rightarrow 0^+}n\sum^{k-1}_{m=1}\frac
  {m^n(1-q^k)}{1-q^m}-k^{n+1}=n \sum^{k-1}_{m=1}m^{n}-k^{n+1} \leq
  0,
\end{equation*}
  where the last inequality follows easily by induction on $k$ and
  this completes the proof.
\end{proof}
\section{Inequalities for the Volume of the Unit Ball in $R^n$}
\label{sec 5} \setcounter{equation}{0}
   In this section, we apply some of our results in the previous sections to study
    inequalities for the volume $\Omega_n$ of the unit ball in
    $R^n$:
 \begin{equation*}
   \Omega_n = \frac {\pi^{n/2}}{\Gamma(1+n/2)}.
 \end{equation*}
   There exists many inequalities involving $\Omega_n$, we refer the
   reader to \cite{alz1.51} and the references therein. In
   \cite{alz1.51}, Alzer has shown that for $n \geq 2$,
 \begin{equation*}
     \frac {2n+1}{4\pi} < \genfrac(){1pt}{}{\Omega_{n-1}}{\Omega_n}^2 \leq \frac
     {2n+\pi-2}{4\pi}.
 \end{equation*}
     The above results were improved by him
     recently \cite[Theorem 3.8]{alz2} to be:
\begin{equation}
\label{1.3}
     \frac 1{2\pi}\exp \Big (\frac \alpha{n}+\psi(n) \Big ) < \genfrac(){1pt}{}{\Omega_{n-1}}{\Omega_n}^2
     \leq \frac 1{2\pi}\exp \Big (\frac \beta{n}+\psi(n) \Big ),
\end{equation}
      with $\alpha=1$ and $\beta=2(\ln(8/\pi)+\gamma-1)$, where $\gamma$ is the Euler's constant
       We point out here the left-hand side inequality in \eqref{1.3}
   is a consequence of the left-hand side inequality of \eqref{3.10}.
    In fact, as in \cite{alz2}, the left-hand side inequality in
   \eqref{1.3} follows from
\begin{equation*}
   2\ln \Gamma(x+1)-2\ln \Gamma(x+1/2)+ \ln 2 - \psi(2x+1)>0
\end{equation*}
   for $x \geq 0$. The above inequality then follows from the left-hand side of
   \eqref{3.10} with $s=1/2$ and the following relation \cite[p. 259]{Ab&S}:
\begin{equation*}
   \psi(2x)=\frac 1{2}\psi(x)+\frac 1{2}\psi(x+1/2)+\ln 2.
\end{equation*}

   Alzer has also shown \cite[Theorem 3]{alz1.51} that for $n \geq
   1$,
\begin{equation*}
   (1+1/n)^{\alpha} \leq \frac{\Omega^2_n}{\Omega_{n-1}\Omega_{n+1}}
   \leq (1+1/n)^{\beta},
\end{equation*}
   with $\alpha=2-\ln \pi/\ln 2, \beta=1/2$ being best possible.
   We now improve the upper bound above in the following
\begin{theorem}
\label{thm5.1}
   For $n \geq 1$,
\begin{equation*}
   (1+ \frac 1{n+1})^{1/2} \leq \frac{\Omega^2_n}{\Omega_{n-1}\Omega_{n+1}}
   \leq (1+ \frac 1{n+1/2})^{1/2},
\end{equation*}
\end{theorem}
\begin{proof}
    Theorem \ref{thm5} implies that $a(x)=1/g(x; 1+(n-1)/2, 1+n/2, 1+(n-1)/2)$
 is logarithmically completely monotonic
   on $(-1-(n-1)/2, +\infty)$. It follows from this and Lemma
   \ref{thm6} that $a(x)/a(x+1/2)$ is logarithmically completely monotonic
   on $(-1-(n-1)/2, +\infty)$. This further implies for
   $x>-1-(n-1)/2$:
\begin{equation*}
   \frac {a(x)}{a(x+1/2)} \geq \lim_{x \rightarrow +\infty} \frac
   {a(x)}{a(x+1/2)}=1,
\end{equation*}
   where the equality above follows from \eqref{3.05}. On taking
   $x=0$ above, we get the first inequality of the theorem.
   Similarly, the consideration of $b(x)=g(x; 1+(n-1)/2, 1+n/2, 1+n/2-3/4)$
  and $b(x)/b(x+1/2)$ gives the second inequality.
\end{proof}

\section*{Acknowledgement}
\setcounter{equation}{0}
   This work was partially carried out while the author was visiting the American Institute of
   Mathematics in fall 2005 and the Centre de Recherches Math\'ematiques
at the Universit\'e de Montr\'eal in spring 2006.
  We would like to thank both the American Institute of Mathematics and the Centre de Recherches Math\'ematiques
at the Universit\'e de Montr\'eal for their generous support and hospitality.



\end{document}